\documentclass[twoside]{article}
\usepackage{amsfonts}
\usepackage{graphicx,color}
\usepackage{amssymb}
\usepackage{amsmath}

\input{tcilatex}

\begin{document}

\date{}

\centerline{}

\centerline {\Large{\bf On Some q-Analogues of the Natural Transform and}}

\centerline {}

\centerline{\Large{\bf  Further Investigations}}

\centerline{}

\centerline{\bf {S. K. Q. Al-Omari}}

\centerline{}

\centerline{Department of Applied Sciences; Faculty of Engineering
Technology; Al-Balqa' } \centerline{ Applied University; Amman 11134; Jordan}
\centerline{s.k.q.alomari@fet.edu.jo}

\centerline{\bf {}}

\centerline{\bf Abstract}

Some $q$-analogues of classical integral transforms have recently been
investigated by many authors in diverse citations. The $q$-analogues of the
Natural transform are not known nor used. In the present paper, we are
concerned with definitions and investigations of the $q$-theory of the
Natural transform and some applications. We present two types of $q$%
-analogues of the cited transform on given sets and get results of the
nominated analogues for certain class of functions of special type. We
declare here that given results are new and they complement recent known
results related to $q$-Laplace and $q$-Sumudu transforms. Over and above, we
present some supporting examples to illustrate effectiveness of the given
results .\medskip 


\section{Introduction}

The study of $q$-analogues of classical integral transforms has not yet been
developed to a great extent. This partially can be explained by the fact
that one is not very familiar with the $q$-theory and that basic $q$%
-integral transforms do not occur frequently in physics. It by no means our
aim to give in this paper a new general results in the theory of $q$%
-calculus. But we restrict our selves to the necessary theory to give some $%
q $-analogues of an integral transform named as the Natural transform and to
estimate the transform relevant properties. The classical theory of the
Natural transform is closely related to the classical theory of Laplace and
Sumudu transforms, two of the best known of all integral transforms. The
Natural transform for functions of exponential order is defined over the set 
$A,$%
\begin{equation}
A=\left\{ f\left( t\right) \left\vert \exists M,\tau _{1}\text{ and/or }\tau
>0\text{: }\left\vert f\left( t\right) \right\vert <Me^{\frac{\left\vert
t\right\vert }{\tau _{j}}},\text{ }t\in \left( -1\right) ^{j}\times \left[
0,\infty \right) ,\text{ }j=1,2\right. \right\}  \tag{1}
\end{equation}%
by the integral equation%
\begin{equation}
N\left( f\right) \left( u;v\right) =\dint\nolimits_{0}^{\infty }f\left(
ut\right) \exp \left( -vt\right) dt,  \tag{2}
\end{equation}%
where $\func{Re}v>0$ and $u\in \left( -\tau _{1},\tau _{2}\right) ,$ $u$ and 
$v$ being the transform variables.

The Natural transform strictly converges to the Sumudu transform (See $\left[
28,29,30\right] $ ) for $v=1$ and it strictly converges to the Laplace
transform for $u=1.$ Further fundamental properties of this transform and
its application to differential equations are given by $\left[ \text{14}%
\right] $ and $\left[ 15,16,17\right] $, respectively$.$

We organize this paper as follows. In section 2, we recall some definitions
and notations from the $q$-calculus. In Section 3, we specify the $q$%
-analogues of the Natural transform in terms of a series representation. In
Section 4, we apply the first of two analogues of the Natural transform to a
certain class of special functions and pick up some results by considering $%
q $-Laplace and $q$-Sumudu transforms. In Section 5, we figure out some
values of the second representation of the $q$-transform and extend the
resulting theorems to the case of $q$-Laplace and $q$-Sumudu transforms. In
Section 6 and 7, we derive certain results concerning Fox's $H_{q}$-function
and propose some counter examples as an application to the previous theory.

\section{The $q$-Calculus}

\noindent For the convenience of the reader, we provide a summary of
mathematical notations used in this paper. Throughout this paper,
wheresoever it appears , $q$ statisfies the condition that $0<q<1.$

\noindent The $q$-calculus begins with the definition of the $q$-analogue $%
d_{q}f\left( x\right) $ of the differential of the function, i.e.,%
\begin{equation*}
d_{q}f\left( x\right) =f\left( x\right) -f\left( qx\right) ,
\end{equation*}%
and the $q$-analogue of the derivative$,$%
\begin{equation*}
\frac{d_{q}f\left( x\right) }{d_{q}x}=\frac{f\left( x\right) -f\left(
qx\right) }{\left( q-1\right) x}.
\end{equation*}%
\noindent On certain additional requirements, the $q$-analogue may be
unique, but sometimes it is useful to consider several $q$-analogues of the
same object.

\noindent The $q$-analogues of an integer $n$ ($q$-integer), a factorial of $%
n$ ($q$-factorial of $n$)$,$ and the binomial coefficient $\binom{n}{k}$($q$%
-binomial coefficient) are respectively given as%
\begin{equation*}
\left[ n\right] _{q}=\frac{1-q^{n}}{1-q},\text{ }\left( \left[ n\right]
_{q}\right) !=\left\{ 
\begin{array}{ll}
\prod\limits_{1}^{n}\left[ k\right] _{q}\text{ }, & n=1,2,3,... \\ 
1\text{ \ \ \ \ \ \ \ \ }, & \text{ }n=0%
\end{array}%
\right. \text{ and }\left[ 
\begin{array}{c}
n \\ 
k%
\end{array}%
\right] _{q}=\dprod\limits_{1}^{n}\frac{1-q^{n-k+1}}{1-q^{k}}.
\end{equation*}

\noindent Clearly, 
\begin{equation*}
\lim_{q\rightarrow 1}\left[ n\right] _{q}=n,\text{ }\lim_{q\rightarrow
1}\left( \left[ n\right] _{q}\right) !=n!\text{ and }\lim_{q\rightarrow 1}%
\left[ 
\begin{array}{c}
n \\ 
k%
\end{array}%
\right] _{q}=\binom{n}{k}.
\end{equation*}%
If $\alpha \in 
\mathbb{C}
,$ then the $q$-analogue of $\alpha $ is given as $\dfrac{1-q^{\alpha }}{1-q}
$ and, it sometimes makes sense when $\alpha $ is not, $\left[ \infty \right]
_{q}=\dfrac{1}{1-q}.$

\noindent If $n\in 
\mathbb{N}
,$ then the $q$-analogue of $\left( x+a\right) ^{n}$ and the derivative are
respectively given as 
\begin{equation}
\left. \left( x+a\right) _{q}^{n}=\prod\limits_{j=0}^{n-1}\left(
x+q^{j}a\right) \text{\ and }D_{q}\left( x+a\right) _{q}^{n}=\left[ n\right]
\left( x+a\right) _{q}^{n-1},\text{ }\left( x+a\right) _{q}^{0}=1\right\} . 
\tag{3}
\end{equation}%
If $x=1$ and $a=x,$ then the above formula makes sense for $n=\infty ,$
giving%
\begin{equation}
\left( 1+x\right) _{q}^{\infty }=\prod\limits_{0}^{\infty }\left(
1+q^{k}x\right) .  \tag{4}
\end{equation}%
The $q$-Jackson integral from $0$ to $a$ is given by Jackson $\left[ 11%
\right] $ as%
\begin{equation}
\int\nolimits_{0}^{a}f\left( x\right) d_{q}x=\left( 1-q\right)
a\sum_{0}^{\infty }f\left( aq^{k}\right) q^{k},  \tag{5}
\end{equation}%
provided the sum converges absolutely.

\noindent The $q$-Jackson integral in a generic interval $\left[ a,b\right] $
is given by $\left[ 11\right] $ 
\begin{equation}
\int\nolimits_{b}^{a}f\left( x\right) d_{q}x=\int\nolimits_{0}^{b}f\left(
x\right) d_{q}x-\int\nolimits_{0}^{a}f\left( x\right) d_{q}x.  \tag{6}
\end{equation}%
The improper integral is defined as $\left[ 5\right] $%
\begin{equation}
\int\nolimits_{0}^{\frac{\infty }{A}}f\left( x\right) d_{q}x=\left(
1-q\right) \sum_{n\in 
\mathbb{Z}
}\frac{q^{k}}{A}f\left( \frac{q^{k}}{A}\right) .  \tag{7}
\end{equation}%
The $q$-analogues of the gamma function are defined by $\left[ 5\right] $%
\begin{equation*}
\left. 
\begin{array}{l}
\Gamma _{q}\left( \alpha \right) =\dint\nolimits_{0}^{\frac{1}{1-q}%
}x^{\alpha -1}E_{q}\left( q\left( 1-q\right) x\right) d_{q}x,\text{ } \\ 
\\ 
_{q}\Gamma \left( \alpha \right) =K\left( A;\alpha \right)
\dint\nolimits_{0}^{\frac{\infty }{A\left( 1-q\right) }}x^{\alpha
-1}e_{q}\left( -\left( 1-q\right) x\right) d_{q}x,%
\end{array}%
\right.
\end{equation*}%
where $\alpha >0$ and that%
\begin{equation}
\left. K\left( A;t\right) =A^{t-1}\dfrac{\left( -q/A;q\right) _{\infty }}{%
\left( -q^{t}/A;q\right) _{\infty }}\dfrac{\left( -A;q\right) _{\infty }}{%
\left( -Aq^{1-t};q\right) _{\infty }},\right.  \tag{8}
\end{equation}%
\begin{equation}
\left( a;q\right) _{n}=\prod\limits_{0}^{n-1}\left( 1-aq^{k}\right) ,\text{ }%
\left( a;q\right) _{\infty }=\underset{n\rightarrow \infty }{\lim }\left(
a;q\right) _{n}\text{ }.  \tag{9}
\end{equation}%
for all $t\in 
\mathbb{R}
.$

\noindent The useful notations we need here are $\left[ 10\right] $%
\begin{equation}
\Gamma _{q}\left( x\right) =\dfrac{\left( q;q\right) _{\infty }}{\left(
q^{x};q\right) _{\infty }}\left( 1-q\right) ^{1-x}\text{ and }\left(
a;q\right) _{t}=\dfrac{\left( a;q\right) _{\infty }}{\left( aq^{t};q\right)
_{\infty }},\text{ }t\in 
\mathbb{R}
.  \tag{10}
\end{equation}%
The $q$-analogues of the exponential function are given as%
\begin{equation}
\left. 
\begin{array}{l}
E_{q}\left( t\right) =\dsum\limits_{0}^{\infty }\left( -1\right) ^{n}\dfrac{%
q^{\frac{n\left( n-1\right) }{2}}}{\left( q;q\right) _{n}}t^{n}=\left(
t;q\right) _{\infty },\text{ }t\in 
\mathbb{C}
\\ 
\\ 
e_{q}\left( t\right) =\dsum\limits_{0}^{\infty }\dfrac{1}{\left( q;q\right)
_{n}}t^{n}=\dfrac{1}{\left( t,q\right) _{\infty }},\text{ }t<1.%
\end{array}%
\right\}  \tag{12}
\end{equation}

\noindent The $q$-analogues of the hypergeometric function are defined in
two ways as%
\begin{eqnarray}
_{r}\phi _{s}\left[ \left. 
\begin{array}{c}
a_{1},a_{2},...,a_{r} \\ 
b_{1},b_{2},...,b_{s}%
\end{array}%
\right\vert q,z\right] &=&\sum_{0}^{\infty }\frac{\left(
a_{1},a_{2},...,a_{r};q\right) _{n}}{\left( b_{1},b_{2},...,b_{s};q\right)
_{n}}\frac{z^{n}}{\left( q;q\right) _{n}}  \TCItag{13} \\
&&\text{and}  \notag \\
_{m-k}\Phi _{m-1}\left[ \left. 
\begin{array}{c}
a_{1},a_{2},...,a_{m-k} \\ 
b_{1},b_{2},...,b_{m-1}%
\end{array}%
\right\vert q,z\right] &=&\sum_{0}^{\infty }\frac{\left(
a_{1},...,a_{m-k};q\right) _{n}}{\left( b_{1},...,b_{m-1};q\right) _{n}}%
\left[ \left( -1\right) ^{n}q^{\binom{n}{2}}\right] ^{k}\frac{z^{n}}{\left(
q;q\right) _{n}},  \notag \\
&&  \TCItag{14}
\end{eqnarray}%
where $\left( a_{1},a_{2},...,a_{p};q\right) _{n}=\prod\limits_{0}^{p}\left(
a_{k},q\right) _{n}.$

\noindent For $\left\vert z\right\vert <2$, some $q$-analogues of the Bessel
function are defined as%
\begin{eqnarray}
J_{v}^{\left( 1\right) }\left( z;q\right) &=&\frac{\left( q^{v+1};q\right)
_{\infty }}{\left( q;q\right) _{\infty }}\left( \frac{z}{2}\right)
^{2}{}_{2}\Phi _{1}\left[ \left. 
\begin{array}{cc}
0 & 0 \\ 
q^{v+1} & 
\end{array}%
\right\vert q,\frac{-z^{2}}{2}\right] =\left( \frac{z}{2}\right)
^{v}\sum_{n=0}^{\infty }\frac{\left( \frac{-z^{2}}{4}\right) ^{n}}{\left(
q;q\right) _{v+n}}\left( q;q\right) _{n},  \notag \\
&&  \TCItag{15}
\end{eqnarray}%
\begin{equation}
J_{v}^{\left( 2\right) }\left( z;q\right) =\tfrac{\left( q^{v+1};q\right)
_{\infty }}{\left( q;q\right) _{\infty }}\left( \tfrac{z}{2}\right) ^{v}%
\text{ }_{0}\Phi _{1}\left[ q^{v+1}\left\vert q,\frac{-q^{v+1}z^{2}}{4}%
\right. \right] =\left( \tfrac{z}{2}\right) ^{v}\sum_{0}^{\infty }\dfrac{%
q^{n\left( n+v\right) }}{\left( q;q\right) _{v+n}\left( q;q\right) _{n}}%
\left( \frac{-z^{2}}{4}\right) ^{n},  \tag{16}
\end{equation}%
\begin{equation}
J_{v}^{\left( 3\right) }\left( z;q\right) =\frac{\left( q^{v+1};q\right)
_{\infty }}{\left( q;q\right) _{\infty }}z^{v}{}_{1}\Phi _{1}\left[ \left. 
\begin{array}{c}
0 \\ 
q^{v+1}%
\end{array}%
\right\vert q,qz^{2}\right] =z^{v}\sum_{0}^{\infty }\frac{\left( -1\right)
^{n}q^{\frac{n\left( n-1\right) }{2}}}{\left( q;q\right) _{v+n}\left(
q;q\right) _{n}}\text{ }.\text{ \ \ \ \ \ \ \ \ \ \ \ \ \ \ }  \tag{17}
\end{equation}

\section{The $q$-Analogues of the Natural Transform}

Theory and applications of $q$-integral transforms are evolving rapidly over
the recent years. Since Jackson $[11]$ presented a precise definition of
so-called $q$-Jackson integral and developed $q$-calculus in a systematic
way. It was well known that, in the literature, there are two types of $q$%
-analogues of integral transforms studied in detail by many authors in the
recent past such as Abdi $\left[ 2\right] ,$ Hahn $\left[ 18\right] ,$
Purohit and Kalla $\left[ 3\right] ,$ Albayrak $\left[ 25\right] $, U\c{c}ar
and Albayrak $\left[ 4\right] ,$ Albayrak et al. $\left[ 5\right] $ and $%
\left[ 6\right] ,$ Yadav and Purohit $\left[ 7\right] ,$ Fitouhi and
Bettaibi $\left[ 8\right] $ and $\left[ 9\right] $ and many others, to
mention but a few.\smallskip

In this section of this paper we deem it proper to give the definition of
the $q$-analogues of the Natural transform as in the following
definition.\smallskip

\noindent \textsc{Definition 1.}\textit{\ Let} $\hat{A}$ \textit{and} $%
\check{A}$ \textit{be defined by}

$%
\begin{array}{c}
\hat{A}=\left\{ f\left( t\right) \left\vert \exists M,\tau _{1}\text{ 
\textit{and/or} }\tau >0\text{: }\left\vert f\left( t\right) \right\vert
<ME_{q}\left( \frac{\left\vert t\right\vert }{\tau _{j}}\right) ,\text{ }%
t\in \left( -1\right) ^{j}\times \left[ 0,\infty \right) ,\text{ }%
j=1,2\right. \right\} \\ 
\text{\textit{and}} \\ 
\check{A}=\left\{ f\left( t\right) \left\vert \exists M,\tau _{1}\text{ 
\textit{and/or} }\tau >0\text{: }\left\vert f\left( t\right) \right\vert
<Me_{q}\left( \tfrac{\left\vert t\right\vert }{\tau _{j}}\right) ,\text{ }%
t\in \left( -1\right) ^{j}\times \left[ 0,\infty \right) ,\text{ }%
j=1,2\right. \right\} ,%
\end{array}%
$

\noindent \textit{respectively. Then, we have the following definitions.}

$\left( \mathtt{i}\right) $ \textit{Over the set} $\hat{A},$ \textit{we
define the }$q$\textit{-analogue of the Natural transform of first type as}%
\begin{equation}
N_{q}\left( f\right) \left( u;v\right) =\frac{1}{\left( 1-q\right) u}%
\dint\nolimits_{0}^{\frac{u}{v}}f\left( t\right) E_{q}\left( q\frac{v}{u}%
t\right) d_{q}t,  \tag{18}
\end{equation}%
\textit{provided the integral exists.}

$\left( \mathtt{ii}\right) $ \textit{Over the set} $\check{A},$ \textit{we
define the }$q$\textit{-analogue of the Natural transform of type two as}%
\begin{equation}
_{q}N\left( f\right) \left( u;v\right) =\frac{1}{\left( 1-q\right) }%
\dint\nolimits_{0}^{\infty }f\left( t\right) e_{q}\left( -\frac{v}{u}%
t\right) d_{q}t,  \tag{19}
\end{equation}%
\textit{when the integral exists.\smallskip }

It seems very benificial to us to notice the following relations%
\begin{equation}
\left. 
\begin{array}{c}
N_{q}\left( f\right) \left( 1;v\right) =\left( L_{q}f\right) \left( v\right)
,\text{ \ \ }_{q}N\left( f\right) \left( 1;v\right) =\left( _{q}Lf\right)
\left( v\right) , \\ 
\\ 
N_{q}\left( f\right) \left( u;1\right) =\left( S_{q}f\right) \left( u\right)
,\text{ \ \ }_{q}N\left( f\right) \left( u;1\right) =\left( _{q}Sf\right)
\left( u\right) ,%
\end{array}%
\right.  \tag{20}
\end{equation}%
where $L_{q}\left( S_{q}\right) $ and $_{q}L\left( _{q}S\right) $ are
respectively the $q$-analogues of the Laplace $\left( \text{Sumudu}\right) $
transforms of first $\left( \text{second }\right) $ types; see, for example, 
$\left[ \text{4}\right] \left( \text{ resp.},\left[ 6\right] \right) $.

\noindent In terms of Jackson integral series representation, the $q$%
-analogue of $\left( 18\right) $ can be expressed as%
\begin{equation}
N_{q}\left( f\right) \left( u;v\right) =\frac{\left( q;q\right) _{\infty }}{v%
}\sum_{k\geq 0}\frac{q^{k}f\left( \dfrac{u}{v}q^{k}\right) }{\left(
q;q\right) _{k}},  \tag{21}
\end{equation}%
whereas, the $q$-analogue of $\left( 19\right) $ can similarly be performed
in terms of that series as 
\begin{equation*}
_{q}N\left( f\right) \left( u;v\right) =\sum_{k\in 
\mathbb{Z}
}^{\infty }\frac{\left( q;q\right) _{\infty }}{v}\frac{f\left( q^{k}\right) 
}{\left( -\dfrac{u}{v}q^{k};q\right) _{\infty }}.
\end{equation*}%
Hence, on parity of the fact $\left( a;q\right) _{k}=\dfrac{\left(
a;q\right) _{\infty }}{\left( aq^{k};q\right) _{\infty }}$ $($ for $a=\dfrac{%
-v}{u}),$ the previous equation has the series representation%
\begin{equation}
_{q}N\left( f\right) \left( u;v\right) =\frac{1}{\left( -\dfrac{v}{u}%
;q\right) _{\infty }}\sum_{k\in 
\mathbb{Z}
}^{\infty }q^{k}f\left( q^{k}\right) \left( -\dfrac{u}{v};q\right) _{k} 
\tag{22}
\end{equation}%
that we shall use in later investigations.

\section{$N_{q}$ of Some Special Functions}

In this section of this paper, we apply the analogue $N_{q}$ for certain
functions of special type and extend the work to $q$-Laplace and $q$-Sumudu
transforms. We assume the functions, unless otherwise stated, are of power
series form,%
\begin{equation}
f\left( x\right) =\sum\limits_{n\geq 0}A_{n}x^{n},  \tag{23}
\end{equation}%
where $A_{n}$ is some bounded sequence .\smallskip

In what follows, we establish the following three main theorems of this
section.\smallskip

\noindent \textsc{Theorem 1.} \textit{Let }$\alpha $\textit{\ be a positive
real number and }$f\left( x\right) =\sum\limits_{n\geq 0}A_{n}x^{n}$\textit{%
. Then, we have}%
\begin{equation}
N_{q}\left( x^{\alpha -1}f\left( x\right) \right) \left( u;v\right) =\frac{%
\left( 1-q\right) ^{\alpha -1}}{v^{\alpha }}u^{\alpha -1}\sum\limits_{n\geq
0}A_{n}\frac{u^{n}}{v^{n}}\left( 1-q\right) ^{n}\Gamma _{q}\left( \alpha
+n\right) .  \tag{24}
\end{equation}%
\noindent \textsc{Proof}\textit{\ Let} $\alpha $ \textit{be a positive real
number and} $f\left( x\right) =\sum\limits_{n\geq 0}A_{n}x^{n}$ \textit{be
given. Then, on aid of} $\left( 21\right) ,$ \textit{we write} 
\begin{eqnarray}
N_{q}\left( x^{\alpha -1}f\left( x\right) \right) \left( u;v\right)  &=&%
\frac{\left( q;q\right) _{\infty }}{v}\sum\limits_{k\geq 0}q^{k}\frac{\left( 
\dfrac{u}{v}q^{k}\right) ^{\alpha -1}f\left( \dfrac{u}{v}q^{k}\right) }{%
\left( q;q\right) _{k}}  \notag \\
&=&\frac{u^{\alpha -1}}{v^{\alpha }}\left( q;q\right) _{\infty
}\sum\limits_{k\geq 0}\frac{q^{\alpha k}}{\left( q;q\right) _{k}}%
\sum\limits_{n\geq 0}A_{n}\left( \dfrac{u}{v}q^{k}\right) ^{n}  \notag \\
&=&\frac{u^{\alpha -1}}{v^{\alpha }}\left( q;q\right) _{\infty
}\sum\limits_{n\geq 0}A_{n}\left( \dfrac{u}{v}\right) ^{n}\sum\limits_{k\geq
0}\frac{q^{k\left( \alpha +n\right) }}{\left( q;q\right) _{k}}.  \TCItag{25}
\end{eqnarray}%
Hence, taking into account $\left( 11\right) ,$ $\left( 25\right) $ simply
reveals%
\begin{eqnarray}
N_{q}\left( x^{\alpha -1}f\left( x\right) \right) \left( u;v\right)  &=&%
\frac{u^{\alpha -1}}{v^{\alpha }}\left( q;q\right) _{\infty
}\sum\limits_{n\geq 0}A_{n}\left( \dfrac{u}{v}\right) ^{n}e_{q}\left(
q^{\alpha +n}\right)   \notag \\
&=&\frac{u^{\alpha -1}}{v^{\alpha }}\left( q;q\right) _{\infty
}\sum\limits_{n\geq 0}A_{n}\dfrac{u^{n}}{v^{n}}\frac{1}{\left( q^{\left(
\alpha +n\right) };q\right) _{\infty }}.  \TCItag{26}
\end{eqnarray}%
Therefore, by aid of the Equations $\left( 10\right) $ and $\left( 26\right) 
$ we get that%
\begin{eqnarray*}
N_{q}\left( x^{\alpha -1}f\left( x\right) \right) \left( u;v\right)  &=&%
\frac{u^{\alpha -1}}{v^{\alpha }}\left( q;q\right) _{\infty
}\sum\limits_{n\geq 0}A_{n}\left( \frac{u}{v}\right) ^{n}e_{q}\left(
q^{\alpha +n}\right)  \\
&=&\frac{\left( 1-q\right) ^{\alpha -1}}{v^{\alpha }}u^{\alpha
-1}\sum\limits_{n\geq 0}A_{n}\frac{u^{n}}{v^{n}}\frac{\Gamma _{q}\left(
\alpha +n\right) }{\left( 1-q\right) ^{-n}} \\
&=&\dfrac{\left( 1-q\right) ^{\alpha -1}}{v^{\alpha }}u^{\alpha
-1}\sum\limits_{n\geq 0}A_{n}\dfrac{u^{n}\left( 1-q\right) ^{n}}{v^{n}}%
\Gamma _{q}\left( \alpha +n\right) .
\end{eqnarray*}%
\noindent This completes the proof of the theorems.\smallskip 

\noindent \textsc{Theorem 2.} \textit{Let }$\alpha $\textit{\ be a positive
real number. Then, the following hold.}

$\left( \mathtt{i}\right) $ \textit{Let} $\Gamma _{q}\left( \alpha \right) $ 
\textit{be the q-gamma function of the first type}$.$\textit{\ Then, we have}%
\begin{equation*}
N_{q}\left( x^{\alpha -1}f\left( x\right) \right) \left( u;v\right) =\dfrac{%
\left( 1-q\right) ^{\alpha -1}}{v^{\alpha }}u^{\alpha -1}\Gamma _{q}\left(
\alpha \right) .
\end{equation*}

$\left( \mathtt{ii}\right) $ \textit{Let} $a\in 
\mathbb{R}
$ \textit{and} $f\left( x\right) =_{m-k}\Phi _{m-1}\left[ \left. 
\begin{array}{c}
a_{1},a_{2},...,a_{m-k} \\ 
b_{1},b_{2},...,b_{m-1}%
\end{array}%
\right\vert q,ax\right] .$ \textit{Then, we have}%
\begin{eqnarray}
\left( N_{q}x^{\alpha -1}f\left( x\right) \right) \left( u;v\right) &=&%
\dfrac{\Gamma _{q}\left( \alpha \right) \left( 1-q\right) ^{\alpha
-1}u^{\alpha -1}}{v^{\alpha }}  \notag \\
&&\text{ }_{m-k+1}\Phi _{m}\left[ \left. 
\begin{array}{c}
a_{1},a_{2},...,a_{m-k}q^{\alpha } \\ 
b_{1},b_{2},...,b_{m-1},0%
\end{array}%
\right\vert q,\frac{au}{v}\right] .  \TCItag{27}
\end{eqnarray}%
\noindent \textsc{Proof of }$\left( \mathtt{i}\right) $ By assuming $A_{0}=1$
and $A_{n}=0,$ $\forall n\geq 1,$ it follows from $\left( 23\right) $ that $%
f\left( x\right) =1.$

Hence, the first part obviously follows.

\noindent \textsc{Proof of} $\left( \mathtt{ii}\right) $ Appealing to the $q$%
-analogue $\left( 14\right) ,f\left( x\right) $ can fairly be written as%
\begin{equation*}
f\left( x\right) =\sum\limits_{n\geq 0}\frac{\left(
a_{1},a_{2},...,a_{m-k};q\right) _{n}}{\left(
b_{1},b_{2},...,b_{m-1};q\right) _{n}}\left( \left( -1\right) ^{n}q^{\frac{%
n\left( n-1\right) }{2}}\right) \frac{q^{n}}{\left( q;q\right) _{n}}x^{n}.
\end{equation*}

On setting%
\begin{equation}
A_{n}=\frac{\left( a_{1},a_{2},...,a_{m-k};q\right) _{n}}{\left(
b_{1},b_{2},...,b_{m-1};q\right) _{n}}\left( \left( -1\right) ^{n}q^{\frac{%
n\left( n-1\right) }{2}}\right) ^{k}a^{n},  \tag{28}
\end{equation}%
we get the power series representation of type $\left( 23\right) .$ Hence,
Theorem 1 reveals%
\begin{equation*}
N_{q}\left( x^{\alpha -1}f\left( x\right) \right) \left( u;v\right) =\frac{%
\left( 1-q\right) ^{\alpha -1}}{v^{\alpha }}u^{\alpha -1}\sum\limits_{n\geq
0}A_{n}\frac{u^{n}}{v^{n}}\left( 1-q\right) ^{n}\Gamma _{q}\left( \alpha
+n\right) .
\end{equation*}%
Therefore, on using the identity of the gamma function $\left[ 1\right] $%
\begin{equation}
\Gamma _{q}\left( x+j\right) =\frac{\left( q^{x};q\right) _{j}}{\left(
1-q\right) _{j}}\Gamma _{q}\left( x\right) ,  \tag{29}
\end{equation}%
the previous equation consequently gives%
\begin{eqnarray*}
N_{q}\left( x^{\alpha -1}f\left( x\right) \right) \left( u;v\right) &=&\frac{%
\left( 1-q\right) ^{\alpha -1}}{v^{\alpha }}u^{\alpha -1}\sum\limits_{n\geq
0}A_{n}\frac{u^{n}}{v^{n}}\left( q^{\alpha };q\right) _{n}\text{ }\Gamma
_{q}\left( \alpha \right) \\
&=&\frac{\Gamma _{q}\left( \alpha \right) \left( 1-q\right) ^{\alpha
-1}u^{\alpha -1}}{v^{\alpha }}\sum\limits_{n\geq 0}A_{n}\frac{u^{n}}{v^{n}}%
\left( q^{\alpha };q\right) _{n}\frac{u^{n}}{v^{n}} \\
&=&\frac{\Gamma _{q}\left( \alpha \right) \left( 1-q\right) ^{\alpha
-1}u^{\alpha -1}}{v^{\alpha }}\text{ }_{m-k+1}\Phi _{m}\left[ \left. 
\begin{array}{c}
a_{1},a_{2},...,a_{m-k}q^{\alpha } \\ 
b_{1},b_{2},...,b_{m-1},0%
\end{array}%
\right\vert q,\frac{au}{v}\right] .
\end{eqnarray*}%
\noindent This completes the proof of the theorem.\smallskip

\noindent An inspection of the previous two main theorems leads to the
following list of inclusions:\smallskip

\noindent On setting $k=m=1,$ Theorem $2\left( \mathtt{ii}\right) \left( 
\text{for }\alpha =1\right) $ yields%
\begin{equation}
N_{q}\left( E_{q}\left( ax\right) \right) \left( u;v\right) =\frac{1}{v}%
\text{ }_{1}\Phi _{1}\left[ \left. 
\begin{array}{c}
q \\ 
0%
\end{array}%
\right\vert q,\frac{au}{v}\right] .  \tag{30}
\end{equation}%
Similarly, by inserting $\alpha =1$, Theorem $2\left( \mathtt{i}\right) $
instantly shows%
\begin{equation}
\left( N_{q}\left( 1\right) \right) \left( u;v\right) =\frac{1}{v}.  \tag{31}
\end{equation}%
Moreover, by setting $\alpha =n+1$ in the first part of the theorem$,$ it
yields 
\begin{equation}
N_{q}\left( x^{n}\right) \left( u;v\right) =\frac{\left( 1-q\right) ^{n}}{%
v^{n+1}}u^{n}\left( \left[ n\right] _{q}\right) !\text{ }.  \tag{32}
\end{equation}%
Therefore, $\left( 31\right) $ and $\left( 32\right) $ when designated reveal%
\begin{equation*}
N_{q}\left( 1+x^{2}+...+x^{n}\right) \left( u;v\right) =\sum\limits_{k\geq 0}%
\frac{\left( 1-q\right) ^{k}}{v^{k+1}}u^{k}\left( \left[ k\right]
_{q}\right) !\text{ }.
\end{equation*}%
\noindent We finally establish the main third theorem of this
section.\smallskip

\noindent \textsc{Theorem 3.} \textit{Let }$f\left( x\right) $\textit{\ be
given as} $f\left( x\right) =$ $_{r}\phi _{p}\left[ \left. 
\begin{array}{c}
a_{1},a_{2},...,a_{r} \\ 
b_{1},b_{2},...,b_{p}%
\end{array}%
\right\vert q,ax\right] $ \textit{and }$\alpha >0.$\textit{\ Then, we have}%
\begin{equation}
N_{q}\left( x^{\alpha -1}f\left( x\right) \right) \left( u;v\right) =\frac{%
\Gamma _{q}\left( \alpha \right) \left( 1-q\right) ^{\alpha -1}u^{\alpha -1}%
}{v^{\alpha }}\text{ }_{r+1}\phi _{p}\left[ \left. 
\begin{array}{c}
a_{1},a_{2},...,a_{r},q^{\alpha } \\ 
b_{1},b_{2},...,b_{p}%
\end{array}%
\right\vert q,a\frac{u}{v}\right] .  \tag{33}
\end{equation}

\noindent \textsc{Proof} By taking into account $\left( 13\right) ,$ $%
f\left( x\right) $ is given the series representation%
\begin{equation*}
f\left( x\right) =\sum\limits_{n\geq 0}\frac{\left(
a_{1},a_{2},...,a_{m-k};q\right) _{n}}{\left(
b_{1},b_{2},...,b_{m-1},q\right) _{n}}\frac{a^{n}}{\left( q;q\right) _{n}}%
x^{n}.
\end{equation*}%
On setting $A_{n}=\dfrac{\left( a_{1},a_{2},...,a_{m-k};q\right) _{n}}{%
\left( b_{1},b_{2},...,b_{m-1},q\right) _{n}}\dfrac{a^{n}}{\left( q;q\right)
_{n}}$ and using $\left( 21\right) $ and $\left( 29\right) ,$ we get 
\begin{eqnarray*}
N_{q}\left( x^{\alpha -1}f\left( x\right) \right) \left( u;v\right) &=&\frac{%
\left( 1-q\right) ^{\alpha -1}u^{\alpha -1}}{v^{\alpha }}\sum\limits_{n\geq
0}A_{n}\frac{u^{n}}{v^{n}}\left( 1-q\right) ^{n}\Gamma _{q}\left( \alpha
+n\right) \\
&=&\frac{\Gamma _{q}\left( \alpha \right) \left( 1-q\right) ^{\alpha
-1}u^{\alpha -1}}{v^{\alpha }}\sum\limits_{n\geq 0}A_{n}\left( q^{\alpha
};q\right) _{n}\frac{u^{n}}{v^{n}} \\
&=&\frac{\Gamma _{q}\left( \alpha \right) \left( 1-q\right) ^{\alpha
-1}u^{\alpha -1}}{v^{\alpha }}\text{ }_{r+1}\phi _{\rho }\left[ \left. 
\begin{array}{l}
a_{1},a_{2},...,a_{r},q^{\alpha } \\ 
b_{1},b_{2},...,b_{p}%
\end{array}%
\right\vert q,a\frac{u}{v}\right] .
\end{eqnarray*}%
This completes the proof of the theorem.\smallskip

\noindent By setting $p=0,\alpha =1$ and $r=0$ in $\left( 33\right) $, the
above theorem leads to%
\begin{equation}
N_{q}\left( e_{q}\right) \left( ax\right) \left( u;v\right) =\frac{1}{v}%
\text{ }_{1}\phi _{0}\left[ q;-\left\vert q,a\frac{u}{v}\right. \right] =%
\frac{1}{v}\sum\limits_{n\geq 0}\left( \frac{au}{v}\right) ^{n}=\frac{1}{v-au%
},\text{ }\left\vert au\right\vert <v.  \tag{34}
\end{equation}%
Further, by fixing $v=1$ and taking account of $\left( 20\right) ,$ $\left(
34\right) $ spreads the result to the case of $q$-Sumudu transform giving 
\begin{equation*}
S_{q}\left( e_{q}\right) \left( ax\right) \left( u\right) =\frac{1}{1-au},%
\text{ }\left\vert au\right\vert <1.
\end{equation*}%
Similarly, by fixing $u=1$ and consulting $\left( 20\right) $ yield the
following case of $q$-Laplace transform 
\begin{equation*}
L_{q}\left( e_{q}\left( ax\right) \right) \left( v\right) =L_{q}\left(
e_{q}\left( ax\right) \right) \left( v\right) =\frac{1}{v-a},\text{ }%
\left\vert a\right\vert <v.
\end{equation*}%
From above investigations, we, further, deduce%
\begin{equation}
N_{q}\left( \sin _{q}ax\right) \left( u;v\right) =N_{q}\left( \dfrac{%
e_{q}\left( iax\right) -e_{q}\left( -iax\right) }{2i}\right) \left(
u;v\right) =\frac{au}{v^{2}+a^{2}u^{2}},\text{ }\left\vert au\right\vert <v,
\tag{35}
\end{equation}

and%
\begin{equation}
N_{q}\left( \cos _{q}ax\right) \left( u;v\right) =N_{q}\left( \dfrac{%
e_{q}\left( iax\right) +e_{q}\left( -iax\right) }{2i}\right) \left(
u;v\right) =\frac{v}{v^{2}+a^{2}u^{2}},\text{ }\left\vert au\right\vert <v. 
\tag{36}
\end{equation}%
\noindent Hence, by virtue of $\left( 35\right) $ and $\left( 36\right) $ we
state without proof the following corrollary.\smallskip

\noindent \textsc{Corrollary 4.} \textit{Let }$a$\textit{\ be a real number.
Then, the following hold.}

$\left( \text{i}\right) S_{q}\left( \sin _{q}ax\right) \left( u\right) =%
\dfrac{au}{1+a^{2}u^{2}},$\ $\left\vert au\right\vert <1;$ $\ \ \ \left( 
\text{ii}\right) S_{q}\left( \cos _{q}ax\right) \left( u\right) =\dfrac{1}{%
1+a^{2}u^{2}},$\ $\left\vert au\right\vert <1,$

$\left( \text{iii}\right) L_{q}\left( \sin _{q}ax\right) \left( v\right) =%
\dfrac{a}{v^{2}+a^{2}},$ $\left\vert a\right\vert <v;$ $\ \ \ \ \ \left( 
\text{iv}\right) L_{q}\left( \cos _{q}ax\right) \left( v\right) =\dfrac{v}{%
v^{2}+a^{2}},$ $\left\vert a\right\vert <v.$

\noindent From Theorem 1 we state and prove the following
corrollary.\smallskip

\noindent \textsc{Corrollary 5.} \textit{Let }$a$\textit{\ be a real number
and }$f\left( x\right) =\dsum\limits_{n\geq 0}A_{n}x^{n}$\textit{. Then, we}$%
\ $\textit{have\ }

$N_{q}\left( x^{\alpha -1}J_{2\mu }^{\left( 1\right) }\left( 2\sqrt{ax}%
;q\right) \right) \left( u;v\right) =\dfrac{\left( 1-q\right) ^{\alpha
-1}u^{\alpha -1}}{v^{\alpha }}\Gamma _{q}\left( \alpha \right)
\dsum\limits_{n\geq 0}A_{n}\dfrac{u^{n}}{v^{n}}\left( q^{\alpha };q\right)
_{n}.$

\noindent \textsc{Proof} Let $a$ be a real number, then by aid of $\left(
15\right) ,$ we consider to write 
\begin{equation*}
J_{2\mu }^{\left( 1\right) }\left( 2\sqrt{ax};q\right) =x^{\mu
}\sum\limits_{n\geq 0}\frac{\left( -1\right) ^{n}a^{\mu +n}}{\left(
q;q\right) _{2\mu +n}\left( q;q\right) _{n}}x^{n}.
\end{equation*}%
By replacing $\alpha $ by $\alpha -\mu -1,$ and setting $A_{n}=\dfrac{\left(
-1\right) ^{n}a^{\mu +n}}{\left( q;q\right) _{2\mu +n}\left( q;q\right) _{n}}%
,$ we, partially, get%
\begin{equation*}
x^{\alpha -1}J_{2\mu }^{\left( 1\right) }\left( 2\sqrt{ax};q\right)
=x^{\alpha -1}\sum\limits_{n\geq 0}A_{n}x^{n}=x^{\alpha -1}f\left( x\right) .
\end{equation*}%
Hence, by Theorem 1, we obtain%
\begin{eqnarray*}
N_{q}\left( x^{\alpha -1}J_{2\mu }^{\left( 1\right) }\left( 2\sqrt{ax}%
;q\right) \right) \left( u;v\right) &=&\frac{\left( 1-q\right) ^{\alpha
-1}u^{\alpha -1}}{v^{\alpha }}\sum\limits_{n\geq 0}A_{n}\frac{u^{n}}{v^{n}}%
\left( q^{\alpha };q\right) _{n}\Gamma _{q}\left( \alpha \right) \\
&=&\frac{\left( 1-q\right) ^{\alpha -1}u^{\alpha -1}}{v^{\alpha }}\Gamma
_{q}\left( \alpha \right) \sum\limits_{n\geq 0}A_{n}\frac{u^{n}}{v^{n}}%
\left( q^{\alpha };q\right) _{n}.
\end{eqnarray*}%
\noindent This completes the proof of the corrollary.\smallskip

\noindent In view of Corrollary 5, we have the following easy
statement.\smallskip

\noindent \textsc{Corrollary 6.} \textit{Let }$a$\textit{\ be a real number.
Then, we\ have}

$N_{q}\left( J_{2\mu }^{\left( 1\right) }\left( 2\sqrt{ax};q\right) \right)
\left( u;v\right) =\dsum\limits_{n\geq 0}A_{n}u^{n}\left( q^{\alpha
};q\right) _{n}.$

\noindent Further, Corrollary 5 is expressed in terms of $q$-Sumudu and $q$%
-Laplace transforms as in the following results.\smallskip

\noindent \textsc{Corrollary}\textbf{\textbf{\ }}\textsc{7}. \textit{Let }$a$%
\textit{\ be a positive real number. Then, the following hold.}

$\left( \text{i}\right) $ $S_{q}\left( x^{\alpha -1}J_{2\mu }^{\left(
1\right) }\left( 2\sqrt{ax};q\right) \right) \left( u\right) =\left(
1-q\right) ^{\alpha -1}u^{\alpha -1}\Gamma _{q}\left( \alpha \right)
\dsum\limits_{n\geq 0}A_{n}u^{n}\left( q^{\alpha };q\right) _{n},$

$\left( \text{ii}\right) $ $L_{q}\left( x^{\alpha -1}J_{2\mu }^{\left(
1\right) }\left( 2\sqrt{ax};q\right) \right) \left( v\right) =\dfrac{\left(
1-q\right) ^{\alpha -1}}{v^{\alpha }}\Gamma _{q}\left( \alpha \right)
\dsum\limits_{n\geq 0}A_{n}\dfrac{1}{v^{n}}\left( q^{\alpha };q\right) _{n},$

\noindent \textit{where }$A_{n}=\left( -1\right) ^{n}\dfrac{a^{\mu +n}}{%
\left( q;q\right) _{2\mu +n}\left( q;q\right) _{n}}.$

\noindent Corrollary 6 extends the results to the case of $q$-Sumudu and $q$%
-Laplace transforms as follows.\smallskip

\noindent \textsc{Corrollary 8}. \textit{Let }$a$\textit{\ be a real number.
Then, the following hold.}

\noindent $\left( \text{i}\right) S_{q}\left( J_{2\mu }^{\left( 1\right)
}\left( 2\sqrt{ax};q\right) \right) \left( u\right) =\sum\limits_{n\geq
0}A_{n}u^{n}\left( q^{\alpha };q\right) _{n};\ \left( \text{ii}\right)
L_{q}\left( J_{2\mu }^{\left( 1\right) }\left( 2\sqrt{ax};q\right) \right)
\left( v\right) =\dfrac{1}{v}\sum\limits_{n\geq 0}A_{n}\dfrac{\left(
q^{\alpha };q\right) _{n}}{v^{n}}.$

\section{$_{q}N$ Transform of Special Functions}

In this section of this paper, we are concerned with the study of the $q$%
-analogue $_{q}N$ of some special functions. We are precisely concerned with
the series representation of the transform and getting some results related
to $q$-Laplace and $q$-Sumudu transforms.\medskip 

\noindent \textsc{Theorem 9.} \textit{Let }$f\left( x\right)
=\sum\limits_{n\geq 0}A_{n}x^{n}$\textit{\ and }$\alpha >0.$\textit{\ Then,
we have}%
\begin{equation}
_{q}N\left( x^{\alpha -1}f\left( x\right) \right) \left( u;v\right) =\left( 
\frac{u}{v}\right) ^{\alpha }\left( 1-q\right) ^{\alpha -1}\Gamma _{q}\left(
\alpha \right) \sum\limits_{n\geq 0}A_{n}\frac{\left( q^{\alpha };q\right)
_{n}}{k\left( \dfrac{u}{v};\alpha +n\right) }\left( \frac{u}{v}\right) ^{n}.
\tag{37}
\end{equation}%
\noindent \textsc{Proof} Let the hypothesis of the theorem be satisfied for
some $\alpha >0.$ Then, on account of $\left( 22\right) ,$ we declare that%
\begin{eqnarray}
_{q}N\left( x^{\alpha -1}f\left( x\right) \right) \left( u;v\right) &=&\frac{%
1}{\left( -\frac{v}{u};q\right) _{\infty }}\sum_{k\in 
\mathbb{Z}
}q^{\alpha k}\sum\limits_{n\geq 0}A_{n}q^{kn}\left( -\frac{v}{u};q\right)
_{k}  \notag \\
&=&\frac{1}{\left( -\frac{v}{u};q\right) _{\infty }}\sum\limits_{n\geq
0}A_{n}\left( \frac{u}{v}\right) ^{\alpha +n}\sum_{k\in 
\mathbb{Z}
}\left( \frac{q^{k}}{\frac{u}{v}}\right) ^{\alpha +n}\left( \frac{-v}{u}%
;q\right) _{k}\text{ }.  \TCItag{38}
\end{eqnarray}%
Hence, by taking into account the fact that%
\begin{equation*}
\Gamma _{q}\left( \alpha \right) =\frac{K\left( A;\alpha \right) }{\left(
1-q\right) ^{\alpha -1}\left( -\left( 1/A\right) ;q\right) _{\infty }}%
\sum_{k\in 
\mathbb{Z}
}\left( \frac{q^{k}}{A}\right) ^{\alpha }\left( \frac{-1}{A};q\right) _{k}
\end{equation*}%
where $K\left( A;\alpha \right) =A^{\alpha -1}\dfrac{\left( -q/\alpha
;q\right) _{\infty }}{\left( -q^{t}/\alpha ;q\right) _{\infty }}\dfrac{%
\left( -\alpha ;q\right) _{\infty }}{\left( -\alpha q^{1-t}/;q\right)
_{\infty }}$ [25, Equ.$\left( 24\right) ]$ (for $A=\dfrac{u}{v}$ and $\alpha
=\alpha +n),$ we get%
\begin{eqnarray*}
_{q}N\left( x^{\alpha -1}f\left( x\right) \right) \left( u;v\right)
&=&\sum_{k\in 
\mathbb{Z}
}A_{n}\left( \frac{u}{v}\right) ^{\alpha +n}\frac{\Gamma _{q}\left( \alpha
+n\right) \left( 1-q\right) ^{\alpha +n-1}}{K\left( \frac{u}{v};\alpha
+n\right) } \\
&=&\left( \frac{u}{v}\right) ^{\alpha }\left( 1-q\right) ^{\alpha -1}\Gamma
_{q}\left( \alpha \right) \sum\limits_{n\geq 0}A_{n}\frac{\left( q^{\alpha
};q\right) _{n}}{K\left( \frac{u}{v};\alpha +n\right) }\left( \frac{u}{v}%
\right) ^{n}.
\end{eqnarray*}%
\noindent This completes the proof of the theorem.

As a straightforward corrollary of Theorem 9, the previous theorem (for $%
A_{0}=1,A_{n}=0,$ for $n\geq 1)$ gives 
\begin{equation}
_{q}N\left( x^{\alpha -1}\right) \left( u;v\right) =\dfrac{\left( \frac{u}{v}%
\right) ^{\alpha -1}\left( 1-q\right) ^{\alpha -1}\Gamma _{q}\left( \alpha
\right) }{K\left( \frac{u}{v};\alpha \right) }.  \tag{39}
\end{equation}

\noindent Also, Equ. $\left( 39\right) $ reveals :\smallskip

\noindent \textsc{Corrollary 10.} \textit{The following hold true.}%
\begin{eqnarray}
\left( \mathtt{i}\right) \text{ \ }_{q}L\left( x^{\alpha -1}\right) \left(
v\right) &=&\dfrac{1}{v^{\alpha -1}}\frac{\left( 1-q\right) ^{\alpha
-1}\Gamma _{q}\left( \alpha \right) }{K\left( \frac{1}{v};\alpha \right) }.%
\text{ \ \ \ \ \ \ \ \ \ \ \ \ \ \ \ \ \ \ \ \ \ \ \ \ \ \ \ \ \ \ \ \ } 
\TCItag{40} \\
\left( \mathtt{ii}\right) \text{ }_{q}S\left( x^{\alpha -1}\right) \left(
u\right) &=&\frac{u^{\alpha -1}\left( 1-q\right) ^{\alpha -1}\Gamma
_{q}\left( \alpha \right) }{K\left( u;\alpha \right) }.  \TCItag{41}
\end{eqnarray}%
Further, $\left( 39\right) $ and $\left( 41\right) $ jointly lead to the
conclusion $\left( _{q}N\left( 1\right) \right) \left( u;v\right) =\dfrac{1}{%
K\left( \frac{u}{v};\alpha \right) }.$ Therefore, we are directed to the
results%
\begin{equation*}
_{q}L\left( 1\right) \left( v\right) =\dfrac{1}{K\left( v;1\right) }\text{
and }\left( _{q}S\left( 1\right) \right) \left( u\right) =\dfrac{1}{K\left(
u;1\right) }.
\end{equation*}

\noindent \textsc{Theorem 11.} \textit{Let }$a$\textit{\ be a real number
and }$f\left( x\right) =$ $_{m-k}\Phi _{m-1}\left[ \left. 
\begin{array}{l}
a_{1},a_{2},...,a_{m-k},q^{\alpha } \\ 
b_{1},b_{2},...,b_{m-1}%
\end{array}%
\right\vert q,ax\right] .$ \textit{Then, we have}%
\begin{eqnarray*}
_{q}N\left( x^{\alpha -1}f\left( x\right) \right) \left( u;v\right) &=&\frac{%
\left( \frac{u}{v}\right) ^{\alpha }\left( 1-q\right) ^{\alpha -1}\Gamma
_{q}\left( \alpha \right) }{K\left( \frac{u}{v};\alpha \right) } \\
&&\text{ \ \ \ \ \ \ \ \ \ \ \ \ \ \ \ \ \ \ }_{m-k+1}\Phi _{m-1}\left[
\left. 
\begin{array}{l}
a_{1},a_{2},...,a_{m-k},q^{\alpha } \\ 
b_{1},b_{2},...,b_{m-1}%
\end{array}%
\right\vert q,\frac{au}{vq^{\alpha }}\right] .
\end{eqnarray*}%
\noindent \textsc{Proof} Let the hypothesis of the theorem be satisfied. A
charity of $\left( 14\right) $ gives%
\begin{equation*}
f\left( x\right) =\sum\limits_{0}^{\infty }\frac{\left(
a_{1},a_{2},...,a_{m-k};q\right) _{n}}{\left(
b_{1},b_{2},...,b_{m-1},q\right) _{n}}\left( \left( -1\right) ^{n}q^{\frac{%
n\left( n-1\right) }{2}}\right) ^{k}\frac{a^{n}}{\left( q;q\right) _{n}}%
x^{n}.
\end{equation*}%
On setting $A_{n}=\dfrac{\left( a_{1},a_{2},...,a_{m-k};q\right) _{n}}{%
\left( b_{1},b_{2},...,b_{m-1},q\right) _{n}}\left( \left( -1\right) ^{n}q^{%
\frac{n\left( n-1\right) }{2}}\right) ^{k}\dfrac{a^{n}}{\left( q;q\right)
_{n}}$ and using Theorem 9 it implies%
\begin{equation*}
_{q}N\left( x^{\alpha -1}f\left( x\right) \right) \left( u;v\right) =\left( 
\frac{u}{v}\right) ^{\alpha }\left( 1-q\right) ^{\alpha -1}\Gamma _{q}\left(
\alpha \right) \sum\limits_{0}^{\infty }A_{n}\frac{\left( q^{\alpha
};q\right) _{n}}{K\left( \frac{u}{v};\alpha +n\right) }\left( \frac{u}{v}%
\right) ^{n}.
\end{equation*}%
By using the fact that $K\left( A,\alpha \right) =q^{\alpha -1}K\left(
A,\alpha -1\right) ,$ the preceding equation gives%
\begin{eqnarray*}
_{q}N\left( x^{\alpha -1}f\left( x\right) \right) \left( u;v\right)
&=&\left( \frac{u}{v}\right) ^{\alpha }\frac{\left( 1-q\right) ^{\alpha
-1}\Gamma _{q}\left( \alpha \right) }{K\left( \frac{u}{v};\alpha \right) }%
\sum\limits_{0}^{\infty }A_{n}\left( q^{\alpha };q\right) _{n} \\
&&\text{ \ \ \ \ \ \ \ \ \ \ \ \ \ }\left( \left( -1\right) ^{n}q^{\frac{%
n\left( n-1\right) }{2}}\right) ^{-1}\left( \frac{-1}{q^{\alpha }}\right)
^{n}\left( \frac{u}{v}\right) ^{n} \\
&=&\frac{\left( \frac{u}{v}\right) ^{\alpha }\left( 1-q\right) ^{\alpha
-1}\Gamma _{q}\left( \alpha \right) }{K\left( \frac{u}{v};\alpha \right) }%
\text{ } \\
&\text{ \ \ }&\text{ \ \ \ \ \ \ \ \ \ \ \ \ \ }_{m-k+1}\Phi _{m-1}\left[
\left. 
\begin{array}{c}
a_{1},a_{2},...,a_{m-k},q^{\alpha } \\ 
b_{1},b_{2},...,b_{m-1}%
\end{array}%
\right\vert q,\frac{au}{vq^{\alpha }}\right] .
\end{eqnarray*}%
Hence the theorem is completely proved.\smallskip

\noindent \textsc{Corrollary 12.} \textit{Let }$a$\textit{\ be a real number
and }$f\left( x\right) $\textit{\ be defined in terms of the }$q$\textit{%
-hypergeometric function}%
\begin{equation*}
f\left( x\right) =_{m-k+1}\Phi _{m-1}\left[ \left. 
\begin{array}{c}
a_{1},a_{2},...,a_{m-k}, \\ 
b_{1},b_{2},...,b_{m-1}%
\end{array}%
\right\vert q;ax\right] .
\end{equation*}%
\textit{Then, the following identities hold.}

\noindent $\left( \mathtt{i}\right) _{q}L\left( x^{\alpha -1}f\left(
x\right) \right) \left( v\right) =\dfrac{\left( 1-q\right) ^{\alpha
-1}\Gamma _{q}\left( \alpha \right) }{K\left( \frac{1}{v};\alpha \right) }\
_{m-k+1}\Phi _{m-1}\left[ \left. 
\begin{array}{l}
a_{1},a_{2},...,a_{m-k},q^{\alpha } \\ 
b_{1},b_{2},...,b_{m-1}%
\end{array}%
\right\vert q,\dfrac{a}{vq^{\alpha }}\right] .$

\noindent $\left( \mathtt{ii}\right) _{q}S\left( x^{\alpha -1}f\left(
x\right) \right) \left( u\right) =\dfrac{u^{\alpha }\left( 1-q\right)
^{\alpha -1}\Gamma _{q}\left( \alpha \right) }{K\left( u;\alpha \right) }$ $%
_{m-k+1}\Phi _{m-1}\left[ \left. 
\begin{array}{l}
a_{1},a_{2},...,a_{m-k},q^{\alpha } \\ 
b_{1},b_{2},...,b_{m-1}%
\end{array}%
\right\vert q,\dfrac{au}{q^{\alpha }}\right] .$

\noindent \textsc{Proof} is straightforward from Theorem 11. Details are
therefore omitted.\smallskip

\noindent Let $\alpha =m=k=1$ and $a>0.$ Then, Corrollary 12 gives%
\begin{equation*}
_{q}N\left( E_{q}\left( ax\right) \right) \left( u;v\right) =\dfrac{u}{%
vK\left( \frac{u}{v},1\right) }\text{ }_{1}\Phi _{0}\left[ 
\begin{array}{c}
q \\ 
\end{array}%
\left\vert q,-\dfrac{au}{qv}\right. \right] .
\end{equation*}%
Hence, it follows%
\begin{equation*}
\left. 
\begin{array}{l}
\left( \mathtt{i}\right) _{q}L\left( E_{q}\left( ax\right) \right) \left(
v\right) =\dfrac{1}{vK\left( \frac{1}{v},1\right) }\text{ }_{1}\Phi _{0}%
\left[ 
\begin{array}{c}
q \\ 
\end{array}%
\left\vert q,-\dfrac{a}{qv}\right. \right] . \\ 
\left( \mathtt{ii}\right) _{q}S\left( E_{q}\left( ax\right) \right) \left(
u\right) =\dfrac{u}{K\left( u,1\right) }\text{ }_{1}\Phi _{0}\left[ 
\begin{array}{c}
q \\ 
\end{array}%
\left\vert q,-\dfrac{au}{q}\right. \right] .%
\end{array}%
\right.
\end{equation*}%
Also, readers may easily verify that%
\begin{equation}
_{1}\Phi _{0}\left[ 
\begin{array}{c}
q \\ 
\end{array}%
\left\vert q,-\frac{au}{qv}\right. \right] =\sum\limits_{0}^{\infty }\left( 
\frac{au}{qv}\right) ^{n}=\dfrac{au}{qv+au},\text{ }\left\vert au\right\vert
<qv.  \tag{42}
\end{equation}%
From above and the fact that $K\left( s,1\right) =1$ we have the following
corrollary.\smallskip

\noindent \textsc{Corrollary 13.} \textit{Let }$a$\textit{\ be a real
number. Then, we have}%
\begin{equation}
_{q}N\left( E_{q}\left( ax\right) \right) \left( u,v\right) =\dfrac{qu}{%
\left( qv+au\right) },\text{ }\left\vert au\right\vert <qv.  \tag{43}
\end{equation}

\noindent Hence, $\left( 43\right) $\ indeed reveals

$\left( \mathtt{i}\right) $ $_{q}L\left( E_{q}\left( ax\right) \right)
\left( v\right) =\dfrac{q}{\left( qv+a\right) },$ $\left\vert a\right\vert
<qv;$ $\ \left( \mathtt{ii}\right) $ $_{q}S\left( E_{q}\left( ax\right)
\right) \left( u\right) =\dfrac{qu}{\left( q+au\right) },$ $\left\vert
au\right\vert <q.$

\noindent It may also be mentioned here that Corrollary 13 and the identities%
\begin{equation}
_{q}\sin x=\frac{E_{q}\left( ix\right) -E_{q}\left( -ix\right) }{2i}\text{
and }_{q}\cos x=\frac{E_{q}\left( ix\right) -E_{q}\left( -ix\right) }{2} 
\tag{44}
\end{equation}%
state,\ without proof, the following result.\smallskip

\noindent \textsc{Corrollary 14.} \textit{Let }$a$\textit{\ be a real
number. Then, we have}

$\left( \mathtt{i}\right) $ $_{q}N\left( _{q}\sin ax\right) \left(
u;v\right) =\dfrac{-qau^{2}}{q^{2}v^{2}+a^{2}u^{2}},$ $\left\vert
au\right\vert <qv.$

$\left( \mathtt{ii}\right) $ $_{q}N\left( _{q}\cos ax\right) \left(
u;v\right) =\dfrac{q^{2}uv}{q^{2}v^{2}+a^{2}u^{2}},$ $\left\vert
au\right\vert <qv.$

\noindent Further, Corrollary 14 suggests to have the following conclusions
proclaimed.

$\left( \mathtt{i}\right) $ $_{q}L\left( _{q}\sin ax\right) \left( v\right) =%
\dfrac{-qa}{q^{2}v^{2}+a^{2}},$ $\left\vert a\right\vert <qv;$ $\ \ \ \ \ \
\left( \mathtt{ii}\right) $ $_{q}L\left( _{q}\cos ax\right) \left( v\right) =%
\dfrac{q^{2}v}{q^{2}v^{2}+a^{2}},$ $\left\vert a\right\vert <qv.$

$\left( \mathtt{iii}\right) $ $_{q}S\left( _{q}\sin ax\right) \left(
u\right) =-\dfrac{qa}{q^{2}+a^{2}u^{2}},$ $\left\vert au\right\vert <qv;$ $%
\left( \mathtt{iv}\right) $ $_{q}S\left( _{q}\cos ax\right) \left( u\right) =%
\dfrac{q^{2}}{q^{2}+a^{2}u^{2}},$ $\left\vert au\right\vert <qv.$

\noindent Further results concerning some other special functions can be
obtained similarly.

\section{$q$-Analogues of $N$-Transforms for the $q$-Fox's $H$-Function}

\noindent Let $\alpha _{j}$ and $\beta _{j}$ be positive integers and $0\leq
m\leq N;$ $0\leq n\leq M.$\ Due to $\left[ 27\right] ,$ the $q$-analogue of
the Fox's $H$-function is given as

$H_{M,N}^{m,n}\left[ x;q\left\vert 
\begin{array}{c}
\left( a_{1},\alpha _{1}\right) ,\left( a_{2},\alpha _{2}\right) ,...,\left(
a_{\mu },\alpha _{M}\right) \\ 
\left( b_{1},\beta _{1}\right) ,\left( b_{2},\beta _{2}\right) ,...,\left(
b_{N},\beta _{N}\right)%
\end{array}%
\right. \right] =$

$\ \ \ \ \ \ \ \ \ \ \ \ \ \ \ \ \ \ \ \ \ \ \ \ \ \ \ \ \dfrac{1}{2\pi i}%
\dint\nolimits_{C}\dfrac{\dprod\limits_{j=1}^{m}G\left( q^{b_{j}-\beta
_{j}s}\right) \dprod\limits_{j=1}^{n}G\left( q^{1-a_{j}+\alpha _{j}s}\right)
\pi x^{s}}{\dprod\limits_{j=m+1}^{N}G\left( q^{1-b_{j}+\beta _{j}s}\right)
\dprod\limits_{j=n+1}^{M}G\left( q^{a_{j}-\alpha _{j}s}\right) G\left(
q^{1-s}\right) \sin \pi s}d_{q}s$%
\begin{equation}
\tag{45}
\end{equation}%
where $G$ is defined in terms of the product%
\begin{equation}
G\left( q^{^{\alpha }}\right) =\dprod\limits_{0}^{\infty }\left( 1-q^{\alpha
-k}\right) ^{-1}=\frac{1}{\left( q^{^{\alpha }},q\right) _{\infty }}. 
\tag{46}
\end{equation}

\noindent The contour $C$ is parallel to $\func{Re}\left( ws\right) =0,$
with indentations in such away all poles of $G\left( q^{b_{j}-\beta
_{j}s}\right) ,$ $1\leq j\leq m,$ are its right and those of $G\left(
q^{1-a_{j}+\alpha _{j}s}\right) ,$ $1\leq j\leq n,$ are the\ left of $C.$
The integral converges if $\func{Re}\left( s\log x-\log \sin \pi s\right)
<0, $ for large values of $\left\vert s\right\vert $ on $C.$ Hence,%
\begin{equation*}
\left\vert \arg \left( x\right) -w_{2}w_{1}^{-1}\log \left\vert x\right\vert
\right\vert <\pi ,\text{ }\left\vert q\right\vert <1,\text{ }\log
q=-w=-w_{1}-iw_{2},
\end{equation*}%
where $w_{1}$ and $w_{2}$\ are real numbers.

\noindent Indeed, for $\alpha _{i}=\beta _{j}=1,$ for all $i,j,$ $\left(
45\right) $ gives the $q$-analogue of the Meijer's $G$-function%
\begin{equation*}
G_{M,N}^{m,n}\left[ x;q\left\vert 
\begin{array}{c}
a_{1},a_{1},...,a_{M} \\ 
b_{1},b_{2},...,b_{N},%
\end{array}%
\right. \right] =\frac{1}{2\pi i}\dint\nolimits_{C}\frac{\dprod%
\limits_{j=1}^{m}G\left( q^{b_{j}-s}\right) \dprod\limits_{j=1}^{n}G\left(
q^{1-a_{j}+s}\right) \pi x^{2}}{\dprod\limits_{j=m+1}^{N}G\left(
q^{1-b_{j}+s}\right) \dprod\limits_{j=n+1}^{M}G\left( q^{a_{j}-s}\right)
G\left( q^{1-s}\right) \sin \pi s}d_{q}s
\end{equation*}%
\noindent where $0\leq m\leq N;$ $0\leq n\leq M$ and $\func{Re}\left( s\log
x-\log \sin \pi s\right) <0$.\smallskip

We have the following main result of this section.\smallskip

\noindent \textsc{Theorem 15.} $\left( \mathtt{i}\right) $ \textit{Let }$%
\lambda $\textit{\ be any complex number and }$k\in \left( 0,\infty \right)
. $\textit{\ The }$q$\textit{-Natural transform }$N_{q}$\textit{\ of the
Fox's }$H_{q}$\textit{-Function is given as}

$N_{q}\left( x^{\lambda }H_{M,N}^{m,n}\left[ \lambda x^{k};q\left\vert 
\begin{array}{l}
\left( a_{1},\alpha _{1}\right) ,\left( a_{1},\alpha _{1}\right) ,...,\left(
a_{M},\alpha _{M}\right) \\ 
\left( b_{1},\beta _{1}\right) ,b_{2},...,\left( b_{N},\beta _{N}\right)%
\end{array}%
\right. \right] \right) \left( u;v\right) =\dfrac{u^{\lambda }}{v^{\lambda
+1}G\left( q\right) }$

$\ \ \ \ \ \ \ \ \ \ \ \ \ \ \ \ \ \ \ \ \ \ \ \ \ \ \ \ \ \ \ \ \ \ \ \ \ \
\ \ \ \ \ \ \ \ \ \ \ \ \ \ \ \ \ \ \ \ \ \ \ \ \ \ \ \ \ \ \ \ \ \ \ \ \ \
\ \ \ \ \ H_{M+1N}^{m,n+1}\left[ \lambda \dfrac{u^{k}}{v^{k}};q\left\vert 
\begin{array}{l}
\left( -\lambda ,k\right) ,...,\left( a_{M},\alpha _{M}\right) \\ 
\left( b_{1},\beta _{1}\right) ,...,\left( b_{N},\beta _{N}\right)%
\end{array}%
\right. \right] .$

$\left( \mathtt{ii}\right) $ \textit{Let }$\lambda $\textit{\ be any complex
number and }$k\in \left( -\infty ,0\right) .$\textit{\ The }$q$\textit{%
-Natural transform }$N_{q}$\textit{\ of the Fox's }$H_{q}$\textit{-Function
is given as}

$N_{q}\left( x^{\lambda }H_{M,N}^{m,n}\left[ \lambda x^{k};q\left\vert 
\begin{array}{l}
\left( a_{1},\alpha _{1}\right) ,\left( a_{1},\alpha _{1}\right) ...,\left(
a_{M},\alpha _{M}\right) \\ 
\left( b_{1},\beta _{1}\right) ,b_{2},...,\left( b_{N},\beta _{N}\right)%
\end{array}%
\right. \right] \right) \left( u;v\right) =\dfrac{u^{\lambda }}{v^{\lambda
+1}G\left( q\right) }$

$\ \ \ \ \ \ \ \ \ \ \ \ \ \ \ \ \ \ \ \ \ \ \ \ \ \ \ \ \ \ \ \ \ \ \ \ \ \
\ \ \ \ \ \ \ \ \ \ \ \ \ \ \ \ \ \ \ \ \ \ \ H_{M,N+1}^{m+1,n}\left[
\lambda \dfrac{u^{k}}{v^{k}};q\left\vert 
\begin{array}{l}
\left( a_{1},\alpha _{1}\right) ,...,\left( a_{M},\alpha _{M}\right) \\ 
\left( 1+\lambda ,-k\right) ,\left( b_{1},\beta _{1}\right) ,...,\left(
b_{N},\beta _{N}\right)%
\end{array}%
\right. \right] .$

\noindent \textsc{Proof} We prove Part $\left( \mathtt{i}\right) $ since
proof of Part $\left( \mathtt{ii}\right) $ is similar. By considering $%
\left( 45\right) ,$ we have

$N_{q}\left( x^{\lambda }H_{M,N}^{m,n}\left[ \lambda x^{k};q\left\vert 
\begin{array}{l}
\left( a_{1},\alpha _{1}\right) ,\left( a_{1},\alpha _{1}\right) ,...,\left(
a_{M},\alpha _{M}\right) \\ 
\left( b_{1},\beta _{1}\right) ,b_{2},...,\left( b_{N},\beta _{N}\right)%
\end{array}%
\right. \right] \right) \left( u;v\right) =$

$\ \ \ \ \ \ \ \ \ \ \ \ \ \ \ \dfrac{1}{2\pi i}\dint\nolimits_{C}\frac{%
\dprod\limits_{j=1}^{m}G\left( q^{b_{j}-\beta _{j}z}\right)
\dprod\limits_{j=1}^{n}G\left( q^{1-a_{j}+\alpha _{j}z}\right) \pi \left(
\lambda \right) ^{z}}{\dprod\limits_{j=m+1}^{N}G\left( q^{1-b_{j}+\beta
_{j}z}\right) \dprod\limits_{j=n+1}^{M}G\left( q^{a_{j}-\alpha _{j}z}\right)
G\left( q^{1-z}\right) \sin \pi z}N_{q}\left( x^{\lambda +kz}\right) \left(
u;v\right) d_{q}z.$%
\begin{equation}
\tag{49}
\end{equation}%
By virtue of Theorem 2, $\left( 49\right) $ gives

$N_{q}\left( x^{\lambda }H_{M,N}^{m,n}\left[ \lambda x^{k};q\left\vert 
\begin{array}{l}
\left( a_{1},\alpha _{1}\right) ,\left( a_{1},\alpha _{1}\right) ,...,\left(
a_{M},\alpha _{M}\right) \\ 
\left( b_{1},\beta _{1}\right) ,b_{2},...,\left( b_{N},\beta _{N}\right)%
\end{array}%
\right. \right] \right) \left( u;v\right) =$

$\ \ \ \ \dfrac{1}{2\pi i}\dint\nolimits_{C}\frac{\dprod\limits_{j=1}^{m}G%
\left( q^{b_{j}-\beta _{j}z}\right) \dprod\limits_{j=1}^{n}G\left(
q^{1-a_{j}+\alpha _{j}z}\right) \pi \left( \lambda \right) ^{z}}{%
\dprod\limits_{j=m+1}^{N}G\left( q^{1-b_{j}+\beta _{j}z}\right)
\dprod\limits_{j=n+1}^{M}G\left( q^{a_{j}-\alpha _{j}z}\right) G\left(
q^{1-z}\right) \sin \pi z}\dfrac{\left( 1-q\right) ^{\lambda +kz}}{%
v^{\lambda +kz+1}}\Gamma _{q}\left( \lambda +kz+1\right) d_{q}z.$

\noindent The fact that $\left( 1-q\right) ^{\lambda +kz}\Gamma _{q}\left(
\lambda +kz+1\right) =\dfrac{G\left( q\right) ^{\lambda +kz}}{G\left(
q\right) }$ gives

$N_{q}\left( x^{\lambda }H_{M,N}^{m,n}\left[ \lambda x^{k};q\left\vert 
\begin{array}{l}
\left( a_{1},\alpha _{1}\right) ,\left( a_{1},\alpha _{1}\right) ,...,\left(
a_{M},\alpha _{M}\right) \\ 
\left( b_{1},\beta _{1}\right) ,b_{2},...,\left( b_{N},\beta _{N}\right)%
\end{array}%
\right. \right] \right) \left( u;v\right) =$

$\ \ \ \ \ \ \ \ \ \ \ \ \ \ \ \ \ \ \dfrac{1}{2\pi i}\dint\nolimits_{C}%
\frac{\dprod\limits_{j=1}^{m}G\left( q^{b_{j}-\beta _{j}z}\right)
\dprod\limits_{j=1}^{n}G\left( q^{1-a_{j}+\alpha _{j}z}\right) \dprod \left(
\lambda \right) ^{z}u^{\lambda +kz}}{\dprod\limits_{j=m+1}^{N}G\left(
q^{1-b_{j}+\beta _{j}z}\right) \dprod\limits_{j=n+1}^{M}G\left(
q^{a_{j}-\alpha _{j}z}\right) G\left( q^{1-z}\right) \sin \pi zv^{\lambda
+kz+1}}\frac{G\left( q^{\lambda +kz}\right) }{G\left( q\right) }d_{q}z$

$\ \ \ \ \ \ \ \ \ \ \ \ \ \ \ \ \ =\dfrac{1}{2\pi i}\frac{u^{\lambda }}{%
v^{\lambda +1}G\left( q\right) }\dint\nolimits_{C}\frac{\dprod%
\limits_{j=1}^{m}G\left( q^{b_{j}-\beta _{j}z}\right)
\dprod\limits_{j=1}^{n}G\left( q^{1-a_{j}+\alpha _{j}z}\right) G\left(
q^{\lambda +kz}\right) \dprod \left( \lambda \frac{u^{k}}{v^{k}}\right) ^{z}%
}{\dprod\limits_{j=m+1}^{N}G\left( q^{1-b_{j}+\beta _{j}z}\right)
\dprod\limits_{j=n+1}^{M}G\left( q^{a_{j}-\alpha _{j}z}\right) G\left(
q^{1-z}\right) \sin \pi z}d_{q}z$

$\ \ \ \ \ \ \ \ \ \ \ \ \ \ \ \ =\dfrac{u^{\lambda }}{v^{\lambda +1}G\left(
q\right) }H_{M+1,N}^{m,n+1}\left[ \lambda \dfrac{u^{k}}{v^{k}};q\left\vert 
\begin{array}{l}
\left( -\lambda _{1},k\right) ,\left( a_{1},\alpha _{1}\right) ,...,\left(
a_{M},\alpha _{M}\right) \\ 
\left( b_{1},\beta _{1}\right) ,...,\left( b_{N},\beta _{N}\right)%
\end{array}%
\right. \right] ,$

\noindent for $k>0$

This completes the proof of Part $\left( \mathtt{i}\right) $ of the theorem.
Proof of Part $\left( \mathtt{ii}\right) $ is quite similar.\medskip

\noindent Hence the theorem is completely proved.\medskip

\noindent \textsc{Corrollary 16.} $\left( \mathtt{a}\right) $ \textit{Let }$%
\lambda $\textit{\ be any complex number and }$k\in \left( 0,\infty \right)
. $\textit{\ Then, we have}

$\left( \mathtt{i}\right) L_{q}\left( x^{\lambda }H_{M,N}^{m,n}\left[
\lambda x^{k};q\left\vert 
\begin{array}{c}
\left( a_{1},\alpha _{1}\right) ,...,\left( a_{M},\alpha _{M}\right) \\ 
\left( b_{1},\beta _{1}\right) ,...,\left( b_{N},\beta _{N}\right)%
\end{array}%
\right. \right] \right) \left( v\right) =\dfrac{1}{v^{\lambda +1}G\left(
q\right) }\times $

$\ \ \ \ \ \ \ \ \ \ \ \ \ \ \ \ \ \ \ \ \ \ \ \ \ \ \ \ \ \ \ \ \ \ \ \ \ \
\ \ \ \ \ \ \ \ \ \ \ \ \ \ \ \ \ \ \ \ H_{M+1,N}^{m,n+1}\left[ \dfrac{%
\lambda }{v^{k}};q\left\vert 
\begin{array}{l}
\left( -\lambda _{1},k\right) ,\left( a_{1},\alpha _{1}\right) ,...,\left(
a_{\mu },\alpha _{M}\right) \\ 
\left( b_{1},\beta _{1}\right) ,...,\left( b_{N},\beta _{N}\right)%
\end{array}%
\right. \right] .$

$\left( \mathtt{ii}\right) S_{q}\left( x^{\lambda }H_{M,N}^{m,n}\left[
\lambda x^{k};q\left\vert 
\begin{array}{c}
\left( a_{1},\alpha _{1}\right) ,...,\left( a_{M},\alpha _{M}\right) \\ 
\left( b_{1},\beta _{1}\right) ,...,\left( b_{N},\beta _{N}\right)%
\end{array}%
\right. \right] \right) \left( u\right) =\dfrac{u^{\lambda }}{G\left(
q\right) }\times $

$\ \ \ \ \ \ \ \ \ \ \ \ \ \ \ \ \ \ \ \ \ \ \ \ \ \ \ \ \ \ \ \ \ \ \ \ \ \
\ \ \ \ \ \ \ \ \ \ \ \ \ \ \ \ \ H_{M+1,N}^{m,n+1}\left[ \lambda
,u^{k};q\left\vert 
\begin{array}{l}
\left( -\lambda _{1},k\right) ,\left( a_{1},\alpha _{1}\right) ,...,\left(
a_{\mu },\alpha _{M}\right) \\ 
\left( b_{1},\beta _{1}\right) ,...,\left( b_{N},\beta _{N}\right)%
\end{array}%
\right. \right] .$

\noindent $\left( \mathtt{b}\right) $ \textit{Let }$\lambda $\textit{\ be
any complex number and }$k\in \left( -\infty ,0\right) .$\textit{\ Then, we
have}

$\left( \mathtt{i}\right) L_{q}\left( x^{\lambda }H_{M,N}^{m,n}\left[
\lambda x^{k};q\left\vert 
\begin{array}{c}
\left( a_{1},\alpha _{1}\right) ,...,\left( a_{M},\alpha _{M}\right) \\ 
\left( b_{1},\beta _{1}\right) ,...,\left( b_{N},\beta _{N}\right)%
\end{array}%
\right. \right] \right) \left( v\right) =\dfrac{1}{v^{\lambda +1}G\left(
q\right) }\times $

$\ \ \ \ \ \ \ \ \ \ \ \ \ \ \ \ \ \ \ \ \ \ \ \ \ \ \ \ \ \ \ \ \ \ \ \ \ \
\ \ \ \ \ \ \ \ \ \ \ \ \ H_{M,N+1}^{m+1,n}\left[ \lambda x^{k};q\left\vert 
\begin{array}{l}
\left( a_{1},\alpha _{1}\right) ,...,\left( a_{M},\alpha _{M}\right) \\ 
\left( 1+\lambda ,-k\right) ,\left( b_{1},\beta _{1}\right) ,...,\left(
b_{N},\beta _{N}\right)%
\end{array}%
\right. \right] .$

$\left( \mathtt{ii}\right) S_{q}\left( x^{\lambda }H_{M,N}^{m,n}\left[
\lambda x^{k};q\left\vert 
\begin{array}{c}
\left( a_{1},\alpha _{1}\right) ,...,\left( a_{M},\alpha _{M}\right) \\ 
\left( b_{1},\beta _{1}\right) ,...,\left( b_{N},\beta _{N}\right)%
\end{array}%
\right. \right] \right) \left( u\right) =\dfrac{u^{\lambda }}{G\left(
q\right) }\times $

$\ \ \ \ \ \ \ \ \ \ \ \ \ \ \ \ \ \ \ \ \ \ \ \ \ \ \ \ \ \ \ \ \ \ \ \ \ \
\ \ \ \ \ \ \ \ \ \ \ \ H_{M,N+1}^{m+1,n}\left[ \lambda u^{k};q\left\vert 
\begin{array}{l}
\left( a_{1},\alpha _{1}\right) ,...,\left( a_{M},\alpha _{M}\right) \\ 
\left( 1+\lambda ,-k\right) ,\left( b_{1},\beta _{1}\right) ,...,\left(
b_{N},\beta _{N}\right)%
\end{array}%
\right. \right] .$

\section{Concrete Examples}

\noindent By virtue of Theorem 15 and the elementary extentions of some $q$%
-analogues of $\sin _{q}x,$ $\cos _{q}x,$ $\sinh _{q}x$ and $\cosh _{q}x$ in
terms of Fox's H-function; see $\left[ 26\right] ,$ we introduce the
following examples.

\noindent \textsc{Example 1.} \textit{Let }$k=2$\textit{\ and }$\lambda =%
\dfrac{\left( 1-q\right) ^{2}}{4}$\textit{\ in Theorem }15\textit{, then we
have}%
\begin{eqnarray*}
N_{q}\left( x^{\dfrac{\left( 1-q\right) ^{2}}{4}}\sin _{q}x\right) \left(
u;v\right) &=&\frac{\sqrt{\pi }\left( 1-q\right) ^{\frac{-1}{2}}G\left(
q\right) u^{\left( 1-q\right) ^{2}/4}}{v^{\frac{\left( 1-q\right) ^{2}}{4}+1}%
} \\
&\text{ \ \ \ \ \ \ \ \ \ \ }&H_{1,3}^{1,1}\left[ \frac{\left( 1-q\right)
^{2}}{4}\frac{u^{2}}{v^{2}};q\left\vert 
\begin{array}{l}
-\frac{\left( 1-q\right) ^{2}}{4},2 \\ 
\left( \frac{1}{2},1\right) ,\left( 0,1\right) ,\left( 1,1\right)%
\end{array}%
\right. \right] .
\end{eqnarray*}

\noindent \textsc{Example 2.} \textit{On setting }$k=2$\textit{\ and }$%
\lambda =\dfrac{\left( 1-q\right) ^{2}}{4}$\textit{\ in Theorem }15\textit{,
we get}%
\begin{eqnarray*}
N_{q}\left( x^{\frac{\left( 1-q\right) ^{2}}{4}}\cos _{q}x\right) \left(
u;v\right) &=&\frac{\sqrt{\pi }\left( 1-q\right) ^{\frac{-1}{2}}G\left(
q\right) u^{\left( 1-q\right) ^{2}/4}}{v^{\frac{\left( 1-q\right) ^{2}}{4}+1}%
} \\
&\text{ \ \ \ \ \ \ \ \ \ \ \ \ \ \ \ \ \ }&H_{1,3}^{1,1}\left[ \frac{\left(
1-q\right) ^{2}}{4}\frac{u^{2}}{v^{2}};q\left\vert 
\begin{array}{l}
\left( -\frac{\left( 1-q\right) ^{2}}{4},2\right) \\ 
\left( 0,1\right) ,\left( \frac{1}{2},1\right) ,\left( 1,1\right)%
\end{array}%
\right. \right] .
\end{eqnarray*}

\noindent \textsc{Example 3. }\textit{On setting }$k=2$\textit{\ and }$%
\lambda =\dfrac{\left( 1-q\right) ^{2}}{4},$\textit{\ in Theorem }15\textit{%
\ we get}%
\begin{eqnarray*}
N_{q}\left( x^{\frac{\left( 1-q\right) ^{2}}{4}}\sinh _{q}x\right) \left(
u;v\right) &=&\frac{\sqrt{\pi }\left( 1-q\right) ^{\frac{-1}{2}}G\left(
q\right) u^{2}}{iv^{2}} \\
&\text{ \ \ \ \ \ \ \ \ \ \ \ \ \ \ \ }&H_{1,3}^{1,1}\left[ \frac{-\left(
1-q\right) ^{2}u^{2}}{4v^{2}};q\left\vert 
\begin{array}{c}
\left( \frac{\left( 1-q\right) ^{2}}{4},2\right) \\ 
\left( \frac{1}{2},1\right) ,\left( 0,1\right) ,\left( 1,1\right)%
\end{array}%
\right. \right] .
\end{eqnarray*}%
\noindent \textsc{Example 4.} \textit{On setting }$k=2,$\textit{\ }$\lambda =%
\dfrac{\left( 1-q\right) ^{2}}{4},$\textit{\ Theorem }15\textit{\ gives}%
\begin{eqnarray*}
N_{q}\left( x^{\frac{\left( 1-q\right) ^{2}}{4}}\cosh _{q}x\right) \left(
u;v\right) &=&\frac{\sqrt{\pi }\left( 1-q\right) ^{\frac{-1}{2}}G\left(
q\right) u^{2}}{v^{2}} \\
&\text{ \ \ \ \ \ \ \ \ \ \ \ \ \ \ }&H_{1,3}^{1,1}\left[ \frac{-\left(
1-q\right) ^{2}u^{2}}{4v^{2}};q\left\vert 
\begin{array}{l}
\left( \frac{\left( 1-q\right) ^{2}}{4},2\right) \\ 
\left( 0,1\right) ,\left( \frac{1}{2},1\right) ,\left( 1,1\right)%
\end{array}%
\right. \right] .
\end{eqnarray*}%
For similar results of $S_{q}$ and $L_{q}$, we set $u=1$ and $v=1$ in the
preceeding Examples.

\end{document}